\documentclass[final,12pt]{elsarticle}

\usepackage{amsmath,amsthm}
\usepackage{amssymb}
\usepackage{amsfonts}
\newcommand{\field}[1]{\mathbb{#1}}
\newtheorem{thm}{Theorem}[section]

\newtheorem{prop}[thm]{Proposition}

\journal{Applied Mathematics Letters}

\begin{document}

\begin{frontmatter}



\title{On the Fourier transform for a symmetric group homogeneous space}


\author{Ramakrishna Kakarala}

\address{School of Computer Engineering\\
Nanyang Technological University\\
Singapore 637798\\
{\tt\small ramakrishna@ntu.edu.sg}\\
}

\begin{abstract}
By using properties of the Young orthogonal representation, this paper derives a simple form for the Fourier transform of permutations acting on the homogeneous space of $n$-dimensional vectors, and shows that the transform requires $2n-2$ multiplications and the same number of additions. 
\end{abstract}


\end{frontmatter}


\section{Introduction}
\label{}
Let $S_n$ denote the symmetric group on $n$ elements, and $S_{n}^{n}$ the subgroup fixing the $n$-th element.  This paper derives a simplification for the Fourier transform of $S_n$ acting on 
$I_{n} = \{1,2,\ldots,n\}$, or equivalently, the coset space $S_n/S_{n}^{n}$.  Fourier analysis of permutations on $I_n$ 
is important for the statistical analysis of ranked data \cite{diaconis}, pattern matching, and other applications. 

To put the aim of the paper in context, it is useful to consider the ordinary Fourier transform.  Let ${\cal H}$ be the $n \times n$ 
unitary matrix with 
entries ${\cal H}_{k,\ell} = (\sqrt{n})^{-1} e^{-j2\pi (k\ell)/n}$. Then $X = {\cal H}x$ is the discrete Fourier transform of the vector $x$.  If $\Delta_d$ is the translation operator that sends $n \mapsto (n+d)\, {\rm mod}\, n$, and $\Phi_d$ 
is the phase shift matrix 
${\rm diag}[1,e^{j2\pi d/n},\ldots,e^{j2\pi d(n-1)/n}]$, then 
\begin{equation}
{\cal H}\Delta_d\, x = \Phi_d X.
\label{eq:phaseshiftfour}
\end{equation}
Similarly, the permutation Fourier transform presented below converts permutations on $I_n$ to group representation ``phase'' shifts.

Fast Fourier transforms on the groups $S_n$ and their homogeneous spaces have been studied previously. In 
particular, by applying the method of Clausen \cite{clausen}, Maslen and Rockmore \cite[Thm 6.5]{maslenrockmore} give an upper bound for the number of operations (either multiplications or additions) on $S_n/S_{n}^{n}$ as $n^3-n^2$.  
Maslen \cite[Thm 3.5]{maslen} improves the bound on the same space to show that, at most, $3n(n-1)/2$ operations are necessary. This paper shows that $2n-2$ operations are sufficient.     

\section{Background for this paper}
\label{sec:background}
We use standard results for permutations \cite{sagan}.  An adjacent transposition is the permutation $\tau_k=(k,k+1)$ that exchanges the $k$-th and $(k+1)$-th elements but leaves all others unchanged.  Every permutation may be written as a product of adjacent transpositions.  

The Fourier transform on $S_n$ relies on the group's irreducible unitary representations,
with ``frequencies'' given by arithmetic partitions.   Let $\nu = (n_1,\ldots,n_q)$ be 
a partition of $n$ with $n_i \geq n_{i+1}$ and
$n_1 + \ldots + n_q = n$; we write $\nu\vdash n$.  For every $\nu\vdash n$ 
there exists an {\it irreducible representation}, denoted $D_{\nu}$. For example, when $\nu = (n)$,
we have $D_{(n)}(\sigma) = 1$ for all $\sigma\in S_n$. For other $\nu$, we use the Young orthogonal representation (YOR)
to construct the matrices. The Fourier transform of $f:S_n\rightarrow\field{C}$ is 
\begin{equation}
F(\nu) = \sum_{\sigma \in S_n} f(\sigma) D_{\nu}(\sigma), \quad \nu \vdash n.
\label{eq:ft}
\end{equation}
For each $\nu$, the coefficient $F(\nu)$ is a $n_{\nu} \times n_{\nu}$ matrix. If $f(\sigma) = g(\delta\sigma)$, i.e., 
$f$ and $g$ are left translates of each other, then, in a manner 
similar to (\ref{eq:phaseshiftfour}), we obtain that $G(\nu) = D_{\nu}(\delta)^t F(\nu)$.  Of particular interest in this paper is the ``fundamental frequency'' of the transform given by the partition $\phi = (n-1,1)$.  The $(n-1)^2$ entries of $D_{\phi}$ are obtained from the YOR as described in detail below.  

It suffices to describe $D_{\phi}$ on the adjacent transpositions $\{\tau_k\}$, for those generate $S_n$.  
Let $D_{\phi}(\tau_1)$ be the $(n-1)$-dimensional matrix  ${\rm diag}[1,1,\ldots,1,-1]$.
For any $m$, let ${\cal I}_{m}$ denote the $m$-dimensional identity matrix,
and for $k=2,\ldots n-1$, let $R_{k}$ be the $2 \times 2$ symmetric matrix
\begin{equation}
R_{k} = 
\left[\begin{array}{cc}
-\frac{1}{k} & \sqrt{1-\frac{1}{k^2}}\\
\sqrt{1-\frac{1}{k^2}} & \frac{1}{k}
\end{array}
\right].
\label{eq:rkk}
\end{equation}
Now, for $k=2,\ldots,n-1$, define $D_{\phi}(\tau_k)$ to be the symmetric, block-diagonal, matrix
\begin{equation}
D_{\phi}(\tau_k) = \left[
\begin{array}{ccc}
{\cal I}_{n-k-1} & 0 & 0\\
0 & R_k &  0\\
0 & 0 & {\cal I}_{k-2}
\end{array}
\right].
\label{eq:defphi}
\end{equation}
It may be verified that the matrices $\{D_{\phi}(\tau_k)\}$ satisfy the Coexeter relations \cite[pg 88]{sagan},
and generate the irreducible YOR for partition $\phi = (n-1,1)$.  Furthermore, note that the 
decomposition of each $\sigma\in S_{n}^{n}$ into 
$\{\tau_k\}$ excludes $\tau_{n-1}$.  Therefore, from (\ref{eq:defphi}), it follows that, with $\oplus$ denoting matrix 
direct sum and $O_{n-2}(\sigma)$ a ($n-2$)-dimensional orthogonal matrix, 
\begin{equation}
D_{\phi}(\sigma) = 1 \oplus O_{n-2}(\sigma),\quad  {\rm for}\,\,\sigma\in S_{n}^{n}.
\label{eq:reduction}
\end{equation}

\section{Fourier analysis on the homogeneous space}

Our goal is to simplify (\ref{eq:ft}) for functions defined on $I_n$.  We may extend each $f$ defined on $I_n$ to 
a corresponding function $\widetilde{f}$ on $S_n$ by 
$\widetilde{f}(\sigma) = f(\sigma(n))$. Note that $\widetilde{f}$ is constant on left cosets of $S_{n}^{n}$ and, therefore,
``band-limited''.  

\begin{prop}
Given any complex-valued function $f$ defined on $I_n$, the Fourier coefficients $\{\widetilde{F}(\nu)\}$ of the 
function $\tilde{f}$ on $S_n$ defined by 
$\widetilde{f}(\sigma) = f(\sigma(n))$ are such that  $\widetilde{F}(\nu) = 0$ unless $\nu = (n)$ or $\nu = \phi = (n-1,1)$. 
\end{prop}

\proof
Since  $\widetilde{f}(\sigma\delta) = \widetilde{f}(\sigma)$ 
for $\delta \in S_{n}^{n}$, we have by (\ref{eq:ft}) that $\widetilde{F}(\nu) = \widetilde{F}(\nu)D_{\nu}(\delta)^t$.  By averaging both sides over $S_{n}^{n}$, we get $\widetilde{F}(\nu) = \widetilde{F}(\nu)Z(\nu)$ where
\begin{equation}
Z(\nu) = \frac{1}{(n-1)!}\sum_{\delta\in S_{n}^{n}} D_{\nu}(\delta)^t.
\label{eq:projection}
\end{equation}
Now, the Branching Rule \cite[Thm 2.8.3]{sagan} shows that for 
 $\nu = \phi = (n-1,1)$ and $\nu = (n)$, the representation $D_{\nu}$ reduces on the subgroup $S_{n}^{n}$ to contain the constant representation, and that no other irreducible representation does so.  By orthogonality, those matrix entries that are not constant on $S_{n}^{n}$ must sum to zero over the subgroup. 
Therefore $Z(\nu) = 0$ if $\nu$ is not $(n)$ or $\phi$. \qed

If $Z(\phi)_{i,j}$ is the $(i,j)$-th element, then from (\ref{eq:reduction}) we have $Z(\phi)_{1,1}=1$, and, by orthogonality, $Z(\phi)_{i,j}=0$ for all other $(i,j)$.  Since $\widetilde{F}(\phi) = \widetilde{F}(\phi)Z(\phi)$ we obtain that $\widetilde{F}(\phi)$ is zero except possibly in the leftmost column.  Hence, the Fourier transform (\ref{eq:ft}) need only be calculated for the partition $(n)$, and for the $n-1$ entries in the left most column of $D_{\phi}$.  Let ${\cal F}$ denote the the linear transformation taking any $n$-dimensional vector $x$ on $I_n$ to its 
$n$ Fourier transform coefficients $\widetilde{X}((n))$, and the leftmost column entries 
$\widetilde{X}(\phi)_{i,1}$ for $i=1$, $2$, $\ldots$, $n-1$. We write 
\begin{equation}
\widetilde{X} = {\cal F}x
\label{eq:matvectransform}
\end{equation} 
to express the transform, now viewed a matrix operation.  The transform (\ref{eq:matvectransform}) requires at most $n^2$ multiplications and $n(n-1)$ additions.  We show below that, in fact, $2n-2$ operations of each kind are sufficient.    

Our result relies on the following $n\times n$ matrix $U$, whose shape is similar to a ``reverse'' upper
Hessenberg matrix: 
\begin{equation}
U = \begin{pmatrix}
      +1 & +1 & +1 & \cdots & +1 & +1\\
      -1 & -1 & -1 & \cdots & -1 & (n-1) \\
       \vdots  & \vdots & \vdots & \vdots & \vdots & \vdots\\
       -1 & -1 & 2 & 0 & \cdots & 0\\
       -1 & +1 & 0 & 0 & \cdots & 0

\end{pmatrix}
\label{eq:easytransformatrix}
\end{equation}
  Define $A = (UU^{t})^{-1/2}$, and let
\begin{equation}
{\cal T} = A U.
\label{eq:orthogonaltransformmatrix}
\end{equation}
It is easily seen that ${\cal T}$ is an orthogonal matrix, and that $A$ is diagonal
with entries $\{\alpha_k\}_{k=1}^{n}$, with $\alpha_1 = A_{1,1}= 1/\sqrt{n}$, and for $k > 1$, we have 
\begin{equation}
\alpha_k = A_{k,k} = \frac{1}{\sqrt{(n-k+1)(n-k+2)}}. 
\label{eq:alpha}
\end{equation}

Let $x$ be any complex-valued $n\times 1$ vector, and let $X = {\cal T} x$.  To each $\sigma \in S_n$, let the $n \times n$ matrix $P(\sigma)$ be the permutation matrix obtained from the identity ${\cal I}_{n}$ with rows permuted by $\sigma$, i.e., $P(\sigma)_{i,j} = {\cal I}_{\sigma(i),j}$.  
Note that $\sigma \mapsto P(\sigma)$ is an antihomomorphism: $P(\sigma\delta) = P(\delta)P(\sigma)$. To see that, note that for any $\alpha$
we have $P(\alpha)e_k = e_{\alpha^{-1}(k)}$ where $e_k$ is $[0,\ldots,0,1,0,\ldots,0]^t$ with $1$ in the $k$-th position.  
If $x$, $y$ are $n\times 1$ vectors, and $y = P(\sigma)x$, then
$y(i)=x(\sigma(i))$ since $P(\sigma)$ is the permutation operator on column vectors.   We now 
establish the following result, comparable to eq.~(\ref{eq:phaseshiftfour}).                                         

\begin{thm}
For every $\sigma\in S_n$ and all $n\times 1$ vectors $x$, we have that 
\[
{\cal T} P(\sigma)\, x =  \left[1 \oplus D_{\phi}(\sigma)^t\right] X
\]
\end{thm}

\proof   We start by proving for any adjacent transposition $\tau_k$ that 
\begin{equation}
{\cal T}P(\tau_k){\cal T}^{t} = 1\oplus D_{\phi}(\tau_k) = 1\oplus D_{\phi}(\tau_k)^t
\label{eq:fortk}
\end{equation}
Note from (\ref{eq:orthogonaltransformmatrix}), (\ref{eq:alpha}), the $m$-th row of ${\cal T}$ for $m>1$ sums to zero, with
the form 
\begin{equation}
[-\alpha_m, -\alpha_m, \ldots, -\alpha_m, (n-(m-1))\alpha_m,0,\ldots,0].
\label{eq:rowm}
\end{equation}
The product ${\cal T}P(\tau_k)$ is the same as ${\cal T}$ but with columns $k$, $k+1$ swapped. By (\ref{eq:rowm}), we 
see that the only rows of ${\cal T} P(\tau_k)$ that are affected by the column swap are as follows: for $k=1$, row $n$ is modified; 
and for $k>1$, rows $n-(k-1)$, $n-(k-2)$ are modified.  Therefore the product ${\cal T}P(\tau_k){\cal T}^{t}$ is the same as the identity ${\cal I}$ in all entries with the following exceptions: when $k=1$, we have that 
$\left[{\cal T}P(\tau_1){\cal T}^{t}\right]_{nn} = -2\alpha_n^2 = -1$; and when $k>1$, we have that the $2\times 2$ submatrix,
whose upper-left corner indices are $(n-(k-1),n-(k-1))$, has the symmetric form 
\begin{equation}
\left[
\begin{array}{cc}
-(k+1)\alpha_{n-(k-1)}^2 & (k^2-1)\alpha_{n-(k-1)}\alpha_{n-(k-2)} \\
(k^2-1)\alpha_{n-(k-1)}\alpha_{n-(k-2)} & (k-1)\alpha_{n-(k-2)}^2
\end{array}
\right]
\end{equation}
Subsituting from (\ref{eq:alpha}), we find that the above simplifies to $R_{k}$ as defined earlier in (\ref{eq:rkk}), thus 
verifying (\ref{eq:fortk}) for $k=1,2,\ldots,n-1$.   

For the general case, note that every $\sigma \in S_n$ may be written as a product of adjacent
transpositions $\sigma = \tau_{k_1}\cdots\tau_{k_m}$. Since $\sigma\mapsto P(\sigma)$ is an anti-homomorphism,
we have that
\begin{equation}
P(\sigma) = P(\tau_{k_1}\cdots\tau_{k_m}) = P(\tau_{k_m})\cdots P(\tau_{k_1}).
\end{equation}
Applying a similarity transformation with ${\cal T}$ yields
\begin{equation}
{\cal T} P(\sigma) {\cal T}^t = {\cal T} P(\tau_{k_m}) {\cal T}^{t} \cdots {\cal T} P(\tau_{k_1}) {\cal T}^{t}.
\end{equation}
On applying (\ref{eq:fortk}) we establish the theorem: 
\begin{equation}
{\cal T} P(\sigma) {\cal T}^t = 
\left[1\oplus D_{\phi}(\tau_{k_m})^t\right] \cdots 
\left[1\oplus D_{\phi}(\tau_{k_1})^t\right] = 1\oplus D_{\phi}(\sigma)^t. 
\end{equation}
\qed 

Note that for the Fourier transform in (\ref{eq:matvectransform}), we also have 
\[
{\cal F}P(\sigma) x = [1\oplus D_{\phi}(\sigma)]^t \widetilde{X}
\]
from the translation property. Since this is true
for all vectors $x$, we must have ${\cal F} = \left[\lambda_1 {\cal I}_{1} \oplus \lambda_{2} {\cal I}_{n-2}\right] {\cal T}$. To see that, note that ${\cal F} = {\cal C} {\cal T}$ for some matrix ${\cal C}$, and, by applying the Theorem above, we see that ${\cal C}$ commutes with all matrices $1\oplus D_\phi$; the result now follows from Schur's lemma \cite[pg 23]{sagan}. 
  
\subsection{Computation of the transform}

The equality ${\cal T} = A U$, combined with the matrix structure in (\ref{eq:easytransformatrix}), simplifies computation. Let 
$a_{x}(n) = x(1)$, $a_{x}(n-1) = x(1) + x(2)$, $\ldots$, $a_{x}(1) = x(1)+x(2)+\cdots+x(n)$. 
Computing all $\{a_{x}(k)\}$ values requires $n-1$ additions due to recursion.   If   
$\hat{X} = U x$ then $\hat{X}(1) = a_{x}(1)$, $\hat{X}(2) = (n-1)*x(n) - a_{x}(2)$,
$\ldots$, $\hat{X}(n) = x(2) - a_{x}(n)$.  Hence, if $a_{x}$ has been computed, computing $\hat{X}$ requires $(n-2)$ multiplies
and $(n-1)$ additions. Now, since $X = {\cal T}\, x = A \hat{X}$, and $A$ is diagonal, we see that computing $X$ from $\hat{X}$ requires an additional $n$ multiplications. In total, computing the transform $X = {\cal T}x$ requires $2n-2$ multiplications and $2n-2$ additions.  Note that computing $\widetilde{X} = {\cal F}x = {\cal C}AUx$ does not require any extra computation as we may premultiply the diagonal matrix ${\cal C}$ with $A$.

\section{Conclusions}

This paper describes a simplification of the Fourier transform on $S_n/S_{n}^{n}$, and shows that the transform requires $2n-2$
multiplications and the same number of additions.



\bibliographystyle{elsarticle-num}
\bibliography{apml_permfourtransform_rev2}

\end{document}